\documentclass[12pt]{article}
\usepackage{graphicx}
\usepackage{color}
\usepackage{amsmath,amsthm,amsfonts,amssymb}
\usepackage{enumerate}
\usepackage{graphicx}
\usepackage{mathrsfs}
\usepackage{authblk}

\theoremstyle{plain}
\newtheorem{thm}{Theorem}[section]

\theoremstyle{definition}

\theoremstyle{remark}
\newtheorem*{rem}{Remark}

\renewcommand{\phi}{\varphi}

\renewcommand{\a}{\alpha}
\renewcommand{\b}{\beta}

\newcommand{\de}{{\rm d}}
\newcommand{\dy}{{\rm d}y}
\newcommand{\dxi}{{\rm d}\xi}

\newcommand{\dt}{{\rm d}t}

\newcommand{\bear}{\begin{eqnarray}}
\newcommand{\enar}{\end{eqnarray}}
\newcommand{\bess}{\begin{eqnarray*}}
\newcommand{\eess}{\end{eqnarray*}}

\newcommand{\bearn}{\begin{eqnarray*}}
\newcommand{\enarn}{\end{eqnarray*}}
\newcommand{\befig}{\begin{figure}}
\newcommand{\eefig}{\end{figure}}
\newcommand{\GL}{Gr\"unwald-Letnikov }
\newcommand{\sgn}{\operatorname{sgn}}

\definecolor{light-gray}{gray}{0.85}

\begin{document}

\title{Finite Difference Methods for the generator of 1D asymmetric alpha-stable
L\'evy motions}

\author[1]{Yanghong Huang\thanks{yanghong.huang@manchester.ac.uk}}
\author[2]{Xiao Wang\thanks{xiaoheda06@163.com}}
\affil[1]{School of Mathematics, The University of Manchester, Manchester, M13
9PL, UK}
\affil[2]{School of Mathematics and Statistics, Henan University
Kaifeng 475001, China}
\date{}


\maketitle

\begin{abstract}
  Several finite difference methods are proposed for the infinitesimal
  generator of 1D asymmetric $\alpha$-stable L\'evy motions,
  based on the fact that the operator becomes a multiplier in the spectral space.
  These methods take the general form of a discrete convolution, and the coefficients
  (or the weights) in the convolution are chosen to approximate the exact multiplier
  after appropriate transform.
  The accuracy and the associated
  advantages/disadvantages are also discussed, providing some guidance on
  the choice of the right scheme for practical problems, like
  in the calculation of mean exit time for random processes governed
  by general asymmetric $\alpha$-stable motions.
\end{abstract}




\section{Introduction}
\label{intro}
While most random dynamical systems studied during the past century
are driven by Gaussian noise, nowadays many complex systems
are modelled with the presence of non-Gaussian noise generated from even
more exotic L\'evy processes~\cite{MR1280932}.

A L\'evy process is a stochastic process with stationary and independent
increments~\cite{MR2512800,Sato-99}. As a generalisation of the
well-known Wiener process (also called Brownian motion), a (scalar)
L\'evy process $\{X_t\}_{t\geq 0}$ is characterised by the L\'evy-Khintchine
theorem, which states that the characteristic function of $X_t$
can be written as $\mathbb{E}\big(\exp(i\xi X_t)|X_0=0
\big)=\exp\big(t\psi(\xi)\big)$. Here, $\psi$ takes the general form
\begin{equation} \label{eq:chpsi}
    \psi(\xi) = i\mu\xi-\frac{1}{2}q^2\xi^2+
    \int_{\mathbb{R}\setminus\{0\}} \big( e^{iy \xi }-1-i y\xi
    \mathbf{1}_{|y|\leq 1} \big)\nu(\de y),
\end{equation}
where $\mu \in \mathbb{R}, q\geq 0$, $\mathbf{1}_S$ is the indicator
function on the set $S$ and $\nu$ is the so called \emph{L\'evy measure}
that satisfies the condition
\[
    \int_{\mathbb{R}\setminus \{0\}} \min(1,y^2)\ \nu(\de y) < \infty.
\]
Intimately related to a L\'evy process is its \emph{infinitesimal
generator} $\mathcal{A}$, defined for {a} smooth function $f$ as
\[
\mathcal{A}f(x)  := \lim_{t\to 0^+} \frac{ \mathbb{E}\big(f(X_t)\mid
X_0=x\big) - f(x)}{t}.
\]
If $X_t$ is a L\'evy process starting from $x$ at time $t=0$, then
$\mathbb{E}\big( \exp(i\xi(X_t-x))|X_0=x\big)=\exp(t\psi(\xi))$,
and the generator $\mathcal{A}f$ can also be represented as~\cite{Duan}
\begin{align} \label{eq:levyinfgen}
  \mathcal{A}f(x)
    = \mu f'(x) + \frac{q^2}{2} f''(x)
    +\int_{\mathbb{R}\setminus\{0\}} \big[
    f(x+y) - f(x) - y \mathbf{1}_{|y|\leq 1} f'(x)\big] \
    \nu(\de y).
\end{align}
If the Fourier transform of $f$ is defined to be
$\hat{f}(\xi) = \mathcal{F}[f](\xi) = \int_{\mathbb{R}} f(x)\exp(-i\xi x)\dxi$, then
 $\mathcal{A}$ becomes a multiplier in the spectral space\footnote{Notice the
minus sign in $\psi(-\xi)$, because of the factor $\exp(-i\xi x)$ used in the definition
of the Fourier transform and $\exp(i\xi x)$ is used in the inverse transform
$f(x) = \frac{1}{2\pi}\int_{\mathbb{R}} \hat{f}(\xi)\exp(i\xi x)\mathrm{d}x$.},
i.e.,
\[
  \widehat{\mathcal{A}f}(\xi) = \psi(-\xi)\hat{f}(\xi).
\]
This operator $\mathcal{A}$ plays a similar role as the adjoint of the classical
Fokker-Plank operator associated with systems driven by Brownian motions,
and appears in the governing non-local differential equations for quantities
that characterise the underlying stochastic processes.
For instance,  the mean exit time $u(x)$ for the expect time of a particle
starting at $x$ and leaving the domain $\Omega$ satisfies the equation
\[
    \mathcal{A} u = -1\qquad  \mbox{ on }\ \Omega,
\]
with the boundary condition $u = 0$ on $\Omega^c$. The adjoint operator
$\mathcal{A}^*$ also arises in other contexts, like the evolution
of the probability density distribution and escape probability~\cite{Duan}.

In the modelling of practical problems driven by non-Gaussian noise,
the final L\'evy process
usually appears as the limit of many small L\'evy processes~\cite{Sato-99}, giving rise
to \emph{$\alpha$-stable L\'evy processes} that enjoy
better scaling properties. In these new processes,
the L\'evy measure is simply
\begin{equation}\label{eq:alphalevymeas}
   \nu({\rm d}y) = \frac{C_+ \chi_{(0,\infty)}(y)
   +C_- \chi_{(-\infty,0)}(y)}{|y|^{1+\alpha}}({\rm d}y),
 \end{equation}
for some non-negative constants $C_+$ and $C_-$, and the stability
parameter $\alpha \in (0,2)$. Using special
integrals~\eqref{eq:intas} in Appendix~\ref{sec:spin}, $\psi(t)$ defined
in~\eqref{eq:chpsi} becomes
\begin{equation}\label{eq:psiCK}
  \psi(t) = i(\mu + K)\xi - \frac{1}{2}q^2\xi^2
  +\Gamma(-\alpha)C_+(-i\xi)^\alpha + \Gamma(-\alpha)C_-(i\xi)^{\alpha}
\end{equation}
where
\[
  K = \begin{cases}
    (C_--C_+)/(1-\alpha), \qquad &\alpha \neq 1,\cr
      (1-\gamma)(C_+-C_-), & \alpha=1,
  \end{cases}
\]
and $\gamma (\approx 0.57721)$ is the Euler-Mascheroni constant. By isolating the real and imaginary 
parts, the function $\psi(t)$ in~\eqref{eq:psiCK} can be written
in the following form commonly seen in the literature
\begin{equation} \label{Charfun}
  \psi(\xi) =
  i(\mu+K)\xi  - \frac{1}{2}q^2\xi^2
  -\sigma^{\a}|\xi|^{\a}\big(1-i\beta \mbox{sgn}(\xi)
    \tan{({\pi \a}/2)}\big),
\end{equation}
with the scaling parameter $\sigma>0$, the skewness parameter $\beta \in
[-1,1]$. The constants in~\eqref{eq:psiCK}
and~\eqref{Charfun} are related by
\begin{equation} \label{C1C2}
     C_+ = C_{\a}\frac{1+\beta}2, \quad
     C_- = C_{\a}\frac{1-\beta}2,\qquad -1 \leq \beta \leq 1,
\end{equation}
with
\begin{equation}  \label{C_alpha}
 C_{\a} =   -\frac{\sigma^\alpha}{\Gamma(-\a)\cos{({\pi \a}/2)}}.
\end{equation}
This connection between~\eqref{eq:alphalevymeas} and~\eqref{eq:psiCK}
turns out to motivate a particular type of finite difference schemes
originated from fractional derivatives~\cite{gorenflo1998random,MR1847209}.

Parallel to the progress made in the theoretical study of L\'evy
processes~\cite{MR2512800}, related fields have been heavily investigated during
the past two decades,  including random simulations of these more exotic
processes~\cite{MR1306279} and  numerical approximations of the underlying
non-local differential equations for quantities like the mean exit time and escape
probability~\cite{Duan}.
Numerous schemes have also been proposed to discretize these more complicated
integro-differential equations, like those originated from random walk
models~\cite{gorenflo1999discrete,gorenflo1998random,MR1847209}
with fractional derivatives~\cite{diethelm2005algorithms},
those based on quadrature of the singular integral
representation~\cite{MR2552219,Ting12,HuangOberman1,METXiao} and harmonic
extension~\cite{MR3348172}.

In this paper, several finite difference schemes for the
infinitesimal generator $\mathcal{A}$ will be
constructed, by exploring their representations in the
appropriate spectral space. For simplicity, only the normalised
operator corresponding to
\[
  \widehat{\mathcal{A}} =\psi(-\xi) =  \begin{cases}-|\xi|^{\alpha}\big(
1+i\beta\sgn(\xi)\tan({\alpha\pi}/{2})\big),\qquad & \alpha \neq 1,\cr
-|\xi|\big(1+i2\beta \sgn(\xi)\ln |\xi| /\pi\big), & \alpha = 1,
\end{cases}
\]
will be considered;
general cases can be adapted easily by adding
discretization of classical first or second order derivatives for convection (if $\mu\neq
0$) or diffusion (if $q\neq 0$), and by applying appropriate scaling (if $\sigma\neq 1$).
To motivate the general form of the schemes,
basic theoretical tools related to semi-discrete Fourier transform
will be reviewed in Section 2.
In Section 3, specific schemes are constructed based on the approximate multipliers
in the spectral space, followed by several numerical experiments in Section 4.
Finally, special functions and integrals, as well as lengthy derivations
of the weights in certain schemes are collected in the appendix.



\section{The discrete operator and semi-discrete Fourier transform}

For local operators like classical integer order derivatives,
Taylor expansion is usually employed to derive difference
schemes and to assess their order of accuracy. However, for nonlocal
operators like fractional order derivatives~\cite{MR0361633} and infinitesimal
generators like $\mathcal{A}$ above, local expansions around the point
of interest soon become inadequate. The right theoretical framework
has to be built upon new tools, like the
\emph{semi-discrete Fourier transform} for functions defined on a uniform infinite lattice.

Let $v$ be a discrete function (usually sampled from a continuous one)
defined on the grid $\mathbb{Z}_h =\{ x_j = hj \mid j \in \mathbb{Z}\}$
with uniform spacing $h$. Then its semi-discrete Fourier transform is defined as
\begin{equation}\label{eq:sdft}
    \hat{v} (\xi) = \mathcal{F}[v](\xi) = h\sum_{j=-\infty}^\infty
    e^{-i\xi x_j} v_j, \qquad \xi \in I_h := \left[ -\frac{\pi}{h},
    \frac{\pi}{h}\right].
\end{equation}
Here the domain is restricted to be the finite interval $I_h$,
as $\hat{v}$ is periodic with period $2\pi/h$.
The inverse transform can also be readily worked out as
\begin{equation}\label{eq:isdft}
    v_j = \mathcal{F}^{-1}[\hat{v}](x_j) =
    \frac{1}{2\pi} \int_{-\pi/h}^{\pi/h} e^{i\xi x_j}
    \hat{v}(\xi)\dxi.
\end{equation}
When the grid size $h$ goes to zero, both~\eqref{eq:sdft} and~\eqref{eq:isdft}
converge to the continuous Fourier transform and inverse Fourier transform, respectively.

Though not as popular as the closely related Fourier transform
and Fourier series, the transform~\eqref{eq:sdft} has been used for a long time
in numerical analysis
like in spectral methods~\cite{MR1776072} and more recently in the discretization
of the generator of symmetric $\alpha$-stable L\'evy processes (or better
known as \emph{Fractional
Laplacian})~\cite{HuangOberman16}. By the well-known
Plancherel type theorem, if $v$ is in
\[
    \mathcal{\ell}^2(\mathbb{Z}_h) = \Big\{
        v: \mathbb{Z}_h \to \mathbb{C}\mid
        h\sum_{j=-\infty}^\infty |v_j|^2 <\infty\Big\},
\]
then its semi-discrete transform belongs to
\[
    L^2(I_h) = \Big\{
        \hat{v}: I_h \to \mathbb{C} \mid
        \int_{-\pi/h}^{\pi/h} |\hat{v}(\xi)|^2 \dxi <\infty
    \Big\}.
\]
Similarly, the following Parseval's identity holds
\[
  \langle u,v\rangle_{\ell^2(\mathbb{Z}_h)} =
  \frac{1}{2\pi}\langle \hat{u},\hat{v}\rangle_{L^2(I_h)},
\]
with the usual inner products:
\[
  \langle u,v\rangle_{\ell^2(\mathbb{Z}_h)}
  =h\sum_{j=-\infty}^\infty \overline{u_j}v_j,\qquad
  \langle \hat{u},\hat{v}\rangle_{L^2(I_h)} =
  \int_{-\pi/h}^{\pi/h} \overline{\hat{u}(\xi)} \hat{v}(\xi)\dxi.
\]
Finally the following convolution theorem is expected~\cite{HuangOberman16}.
\begin{thm}[Convolution Theorem] \label{lem:conv}
  Let $v$ and $w$ be functions in $\ell^2(\mathbb{Z}_h)$
  together with their semi-discrete transforms $\hat{v}$ and $\hat{w}$
  in $L^2(I_h)$. Then the operator $\mathcal{D}: \ell^2(\mathbb{Z}_h)
  \to \ell^2(\mathbb{Z}_h)$ defined by
  \[
    (\mathcal{D}v)_j = h \sum_{k=-\infty}^\infty v_{j-k}w_k
  \]
  becomes a multiplier in the spectral space $L^2(I_h)$, that is,
  \[
    \widehat{\mathcal{D}v}(\xi) = \hat{w}(\xi)\hat{v}(\xi).
  \]
\end{thm}

Once the foundations have been firmly laid,
finite difference schemes for the infinitesimal generator $\mathcal{A}$
can be constructed based on their representations in the spectral space:
since $\mathcal{A}$ corresponds to the multiplier $\psi(-\xi)$, any reasonable scheme is
expected to be a multiplier $\psi_h(-\xi)$ in the spectral space
that approximates $\psi(-\xi)$. Therefore,
if the discrete operator $\mathcal{A}_hu$ is represented as
$\widehat{\mathcal{A}_h u}(\xi)=\psi_h(-\xi)\hat{u}(\xi)$ in the spectral space,
then in the physical space the scheme reads
\begin{equation} \label{eq:fd}
  \mathcal{A}_h u_j = \sum_{k=-\infty}^\infty w_{j-k}u_k,
\end{equation}
for some weights $\{w_j\}_{j=-\infty}^\infty$ with $\hat{w}(\xi)=h\psi_h(-\xi)$.
The equivalence between these two representations
can be established by the Convolution Theorem~\ref{lem:conv},
and schemes can be constructed directly from appropriate multiplier $\psi_h(-\xi)$,
or alternatively from weights with suitable $\psi_h(-\xi)$ that approximates $\psi(-\xi)$.
Below we focus on the case $\alpha \neq 1$, and the case $\alpha=1$ will be
treated separately.

Before discussing explicit weights detailed in the next section, we make
several observations and assumptions to simplify the presentation. Motivated
from $\psi(\xi) = -|\xi|^{\alpha}\big(1-i\beta \sgn(\xi)\tan({\alpha\pi}/{2})\big)$,
we can assume that the most general form of $\psi_h$ reads
\begin{equation*}
  \psi_h(\xi) = -\tilde{M}_{eh}(\xi) + i\beta
  \tilde{M}_{oh}(\xi)\sgn(\xi)\tan({\alpha\pi}/{2}),
\end{equation*}
where both $\tilde{M}_{eh}(\xi)$ and $\tilde{M}_{oh}(\xi)$
behave like $|\xi|^\alpha$ near the origin.
Because of the consistency  in the accuracy and the appearance of related special integrals
in the expressions of the weights, it happens that in many cases
$\tilde{M}_{eh}(\xi)=\tilde{M}_{oh}(\xi)$, but there is no
need to impose further restrictions like this.
Following the above assumption for $\psi_h$,
\begin{equation}\label{eq:whatM}
\hat{w}(\xi) = h\psi_h(-\xi)
= -h\tilde{M}_{eh}(\xi) - ih\beta \tilde{M}_{oh}(\xi)\sgn(\xi)\tan({\alpha\pi}/{2})
\end{equation}
and from the inverse transform~\eqref{eq:isdft}
\begin{align*}
  w_j &= \frac{1}{2{\pi}}\int_{-\pi/h}^{\pi/h}
  e^{i\xi x_j}\Big[
 -h\tilde{M}_{eh}(\xi) - ih\beta
 \tilde{M}_{oh}(\xi)\sgn(\xi)\tan\left(\frac{\alpha\pi}{2}\right)
 \Big] \dxi \cr
 &= -\frac{h^{-\alpha}}{2\pi}\int_{-\pi}^\pi e^{ij\xi} \left[
   h^{\alpha}\tilde{M}_{eh}\left(\frac{\xi}{h}\right)
   + i\beta h^{\alpha} \tilde{M}_{oh}\left(\frac{\xi}{h}\right)
   \sgn(\xi)\tan\left(\frac{\alpha\pi}{2}\right)
 \right] \dxi.
\end{align*}
Next, because $\mathcal{A}$ is a spatial derivative of order $\alpha$,
the only dependence of the weights $\{w_j\}$
on the grid size $h$ is expected to be the factor $h^{-\alpha}$.
This observation implies that the rescaled symbols or multipliers
$M_e(\xi) := h^{\alpha}\tilde{M}_{eh}(\xi/h)$ and
$M_o(\xi) := h^{\alpha}\tilde{M}_{oh}(\xi/h)$
can be assumed to be independent on $h$.
Under this final assumption, the above integral for $w_j$ can be further simplified
by isolating the real and imaginary parts of the integrand. That is,
\begin{align}\label{eq:weight}
  w_j &= -\frac{h^{-\alpha}}{2\pi} \int_{-\pi}^\pi e^{ij\xi} M_e(\xi)\dxi
  - i\frac{h^{-\alpha}}{2\pi}\beta \tan \left(\frac{\alpha\pi}{2} \right)
  \int_{-\pi}^{\pi} e^{ij\xi} M_o(\xi) \sgn (\xi) \dxi \cr
  &= -\frac{h^{-\alpha}}{\pi}\int_0^\pi M_e(\xi)\cos(j\xi) \dxi
  + \frac{h^{-\alpha}}{\pi} \beta \tan \left(\frac{\alpha\pi}{2}\right)
  \int_0^\pi M_o(\xi) \sin (j\xi) \dxi.
\end{align}
Therefore, the weight $w_j$ consists of two Fourier integrals,
the symmetric part $\int_0^\pi M_e(\xi)\cos(j\xi)\dxi$ and
the anti-symmetric part $\int_0^\pi M_o(\xi)\sin(j\xi)\dxi$;
the symmetric part corresponds to the well-known Fractional Laplacian
studied thoroughly in~\cite{HuangOberman16} under the same framework.


\begin{rem} From the conditions $M_e(\xi) \sim |\xi|^\alpha$ and $M_o(\xi)\sim
  |\xi|^{\alpha}$, we have $\tilde{M}_{eh}(0)=\tilde{M}_{oh}(0)=0$ and
  \[
    \sum_{k=-\infty}^\infty w_k
    =\left.\left(\sum_{k=-\infty}^\infty e^{-i\xi x_k} w_k\right)\right|_{\xi=0}
    = \frac{\hat{w}(0)}{h}
    =\left. -\Big(\tilde{M}_{eh}(\xi) + i\beta
    \tilde{M}_{oh}(\xi)\sgn(\xi)\tan({\alpha\pi}/{2})\Big) \right|_{\xi=0} = 0.
  \]
  As a result, the general scheme~\eqref{eq:fd} can also be written as
  \[
    \mathcal{A}_h u_j = \sum_{k=-\infty}^\infty w_{j-k}(u_k-u_j).
  \]
  This is the form adopted in~\cite{HuangOberman16}, for its
  close resemblance with the singular integral representation of the
  Fractional Laplacian operator.
\end{rem}

\section{Finite difference discretization of the generator $\mathcal{A}$}

In this section, several finite difference schemes
will be constructed with appropriate rescaled multipliers $M_e(\xi)$ and
$M_o(\xi)$,
such that the semi-discrete transform of the weights $\{w_j\}$
in the general scheme~\eqref{eq:fd} becomes
\begin{equation*}
  \hat{w}(\xi) = -h^{1-\alpha}M_e(h\xi) -
  ih^{1-\alpha}\beta M_o(h\xi) \sgn(\xi) \tan \left(\frac{\alpha\pi}{2}\right),
\end{equation*}
where the substitutions $\tilde{M}_{eh}(\xi)=h^{-\alpha}M_e(h\xi)$
and $\tilde{M}_{oh}(\xi)=h^{-\alpha}M_o(h\xi)$ are used in~\eqref{eq:whatM}.
Once $M_e(\xi)$ and $M_o(\xi)$ are chosen,
the weights $\{w_j\}$ are then expressed explicitly as in~\eqref{eq:weight},
and vice versa.

Although no restrictions placed on $M_e(\xi)$ and $M_o(\xi)$ other than the behaviour
$|\xi|^\alpha$ near the origin, there is one main constraint
in practice: the Fourier integrals in~\eqref{eq:weight} can not
be evaluated in closed form expressions in general, while numerical quadrature
could introduce large error, especially when the index $j$ of $w_j$ is large.
In this section, we document several choices of $M_e(\xi)$ and $M_o(\xi)$,
mainly motivated from schemes for the symmetric
fractional Laplacian in~\cite{HuangOberman16}. The key features of these schemes
are explicit weights expressed with special functions available in
standard packages like MATLAB, making any numerical quadrature of
highly oscillatory integrals dispensable.



\subsection{Spectral weights}

The most natural choice is $M_{e}(\xi)=M_{o}(\xi) = |\xi|^\alpha$,
such that $\psi_h(\xi)$ coincides with $\psi(\xi)$ on the interval $[-\pi/h,\pi/h]$.
Consequently, the resulting weights are denoted by $\{w_j^{SP}\}$,
called \emph{spectral weights}. The corresponding scheme has been studied
for the symmetric fractional Laplacian ($\beta=0$):
the expressions $\alpha=2$ (the classical second order derivative)
appeared in the context of spectral
methods~\cite{MR1776072}; general cases for $\alpha \in (0,2)$
were extensively discussed in~\cite{HuangOberman16} under a similar framework.

From~\eqref{eq:weight}, the weights $\{w_j^{SP}\}$  depend
essentially on $\int_0^\pi \xi^\alpha \cos (j\xi)\dxi$ and
$\int_0^\pi \xi^{\alpha}\sin(j\xi)\dxi$, which can be simplified using series
expansions of the trigonometric functions. For instance,
\begin{equation*}
  \int_0^\pi \xi^\alpha \cos (j\xi) \dxi
  =\sum_{n=0}^\infty \frac{(-1)^n j^{2n}}{(2n)!} \int_0^\pi
  \xi^{2n+\alpha}\dxi
  =\sum_{n=0}^\infty \frac{(-1)^nj^{2n}}{(2n)!}
  \frac{\pi^{2n+\alpha+1}}{2n+\alpha+1}.
\end{equation*}
The last series can be represented using the
generalised hypergeometric function ${ }_1F_2$ (not the more popular
Gauss hypergeometric function ${ }_2F_1$, see the Appendix~\ref{sec:spint} for more
details), leading to
\[
  \int_0^\pi \xi^\alpha \cos (j\xi) \dxi
    =\frac{\pi^{\alpha+1}}{\alpha+1}
    \sum_{n=0}^\infty \frac{\left(\frac{\alpha+1}{2}\right)_n}
    {n!\left(\frac{1}{2}\right)_n
    \left(\frac{\alpha+3}{2}\right)_n}\left(-\frac{j^2\pi^2}{4}\right)^n
    =\frac{\pi^{\alpha+1}}{\alpha+1}{ }_1F_2\left(
    \frac{\alpha+1}{2};\frac{1}{2},\frac{\alpha+3}{2};-\frac{j^2\pi^2}{4}
    \right),
\]
where $(a)_n=a(a+1)\cdots(a+n-1)=\Gamma(a+n)/\Gamma(a)$
is the Pochhammer symbol. Similarly,
\begin{align*}
  \int_0^\pi \xi^\alpha \sin (j\xi) \dxi &=
  \sum_{n=0}^\infty \frac{(-1)^nj^{2n+1}}{(2n+1)!}
  \frac{\pi^{2n+\alpha+2}}{2n+\alpha+2} \cr
 & = j\frac{\pi^{\alpha+2}}{\alpha+2}
  \sum_{n=0}^\infty \frac{\left(\frac{\alpha}{2}+1\right)_n}
  {n!\left(\frac{3}{2}\right)_n\left(\frac{\alpha}{2}+2\right)_n}
  \left(-\frac{j^2\pi^2}{4}\right)^n \cr
  &= j\frac{\pi^{\alpha+2}}{\alpha+2} { }_1F_2\left(
  \frac{\alpha}{2}+1;\frac{3}{2},\frac{\alpha}{2}+2;
  -\frac{j^2\pi^2}{4}\right).
\end{align*}
Therefore, from these explicit expressions of the two integrals,
\begin{align*}
w_j^{SP} &= -\frac{h^{-\a}}{\pi}\left[\int_0^\pi
\xi^\a\cos(j\xi)\dxi-\beta\tan\left(\frac{\pi\a}{2}\right)
\int_0^\pi\xi^\a\sin(j\xi)\dxi\right]        \cr
& = -\frac{h^{-\alpha}\pi^\alpha}{\alpha+1}
{ }_1F_2\Big(\frac{\alpha+1}{2}; \frac{1}{2},\frac{\alpha+3}{2};
-\frac{j^2\pi^2}{4}\Big) \cr
&\qquad  +\frac{h^{-\alpha}\pi^{\alpha+1}}{\alpha+2}
  j\beta\tan\left(\frac{\alpha\pi}{2}\right) { }_1F_2\Big(
  \frac{\alpha}{2}+1;\frac{3}{2},\frac{\alpha}{2}+2;-\frac{j^2\pi^2}{4}\Big).
\end{align*}

For the special case $\a=1$ is even simpler: $w_{0}^{SP} = -\pi/2h$ and for $j\neq 0$,
\begin{align*}
    w_j^{SP}     &= -\frac{h^{-1}}{\pi} \int_0^\pi\xi\cos(j\xi)\dxi
    +\frac{2h^{-1}\b}{\pi^2}\int_0^\pi\xi \ln\xi \sin(j\xi)\dxi \\
    &= \frac{1-(-1)^j}{\pi j^2 h}
    +\frac{2\b}{j^2\pi^2h}\Big[ (-1)^j j\pi \ln(\pi)+\text{Si}(j\pi)\Big],
\end{align*}
where $\text{Si}(z)=\int_0^z\frac{\sin{t}}{t} \dt$ is the Sine integral.

\subsection{Gr\"unwald-Letnikov weights and fractional derivatives}
Although the above scheme with the exact multiplier $M_e(\xi)=M_o(\xi)=|\xi|^\alpha$
is natural, {the most popular one in the literature
is based on classical \GL finite differences for Riemann-Liouville fractional derivatives
defined on bounded intervals~\cite{diethelm2005algorithms,MR0361633}. For functions
defined on the whole real line, the original \GL differences can easily be generalised
 by allowing integration limits in the definition going to infinity, leading to the so called
Weyl fractional derivatives.}
Here we follow the historical development of this method and quote the weights first,
before showing the multipliers.
For $\alpha \in (0,1)$, the Weyl fractional derivatives $\mathcal{D}_+^\alpha$ and
$\mathcal{D}_-^\alpha$ for a smooth function $f$ on the real line are defined as~\cite{MR1347689}
\begin{align*}
  \mathcal{D}_\pm^{\alpha} f(x) = \pm \frac{1}{\Gamma(1-\alpha)}
  \frac{\de }{\de x}\int_0^\infty t^{-\alpha}f(x\mp t)\dt
  =-\frac{1}{\Gamma(-\alpha)}\int_0^\infty
  \frac{f(x)-f(x\mp t)}{t^{\alpha+1}}\dt.
\end{align*}
Numerous schemes have been proposed to approximate fractional derivatives like these,
including the classical \GL finite differences
(see~\cite{diethelm2005algorithms,MR0361633} for more details):
\begin{equation}
  \label{eq:glda}
  \mathcal{D}_{\pm}^\alpha f(x_j) \approx
  \mathcal{D}_{h\pm}^\alpha f_j
  :=h^{-\alpha}\sum_{k=0}^\infty (-1)^k\binom{\alpha}{k}f_{j\mp k},
\end{equation}
where  $f_j = f(x_j)$ and
$\binom{\alpha}{k} =
\frac{\Gamma(k-\alpha)}{k!\Gamma(-\alpha)} =
(-1)^k\frac{\alpha(\alpha-1)\cdots(\alpha-k+1)}{k!}$
is the generalised binomial coefficients.
The infinitesimal generator\footnote{The term $-iy\xi\mathbf{1}_{|y|\leq 1}$ is omitted in
the expression of the integral, contributing to the constant $iK\xi$ in~\eqref{eq:psiCK}.}
$\mathcal{A}$ is then related to these fractional derivatives by
\begin{align}\label{eq:fraclevy}
  \mathcal{A}f(x) &= \int_{\mathbb{R}\setminus\{0\}}
  [f(x+y)-f(x)]\nu(\dy) \cr
  &= \Gamma(-\alpha)\Big[
    C_+\mathcal{D}_-^\alpha f(x) + C_-\mathcal{D}_+^\alpha f(x)
  \Big] \cr
  &= -\frac{1}{2\cos(\alpha\pi/2)}\Big[
    (1+\beta)\mathcal{D}_-^\alpha f(x) + (1-\beta)\mathcal{D}_+^\alpha f(x)
  \Big],
\end{align}
where the relations~\eqref{C1C2} and~\eqref{C_alpha} are used in the last step.
As a result, we  obtain a scheme for $\mathcal{A}$ with the fractional derivatives
$\mathcal{D}_\pm^\alpha$ approximated by~\eqref{eq:glda},
and the weights $\{w_j^{GL}\}$, called \GL weights, are collected as
\begin{equation*}
  h^{\alpha} w^{GL}_j = \begin{cases}
    -\frac{1+\beta}{2\cos(\alpha\pi/2)}(-1)^j \binom{\alpha}{j},\qquad & j>0,\cr
    -\frac{1}{\cos(\alpha\pi/2)}, & j=0,\cr
    -\frac{1-\beta}{2\cos(\alpha\pi/2)}(-1)^j \binom{\alpha}{-j},
    \qquad & j<0.
  \end{cases}
\end{equation*}
It is easy to check that $w_j^{GL}$ is positive for
all $\alpha \in (0,1)$ and $j\neq 0$.
The associated rescaled multipliers $M_e(\xi)$ and $M_o(\xi)$ can then be obtained
by using the definition of the semi-discrete
transform~\eqref{eq:sdft} and by recognizing that the (positive)
weights are proportional to the binomial coefficients of
$(1-z)^{\alpha}$, i.e.,
\begin{align*}
  M_e(\xi) &=     \frac{1}{2\cos (\alpha \pi/2)}
    \Big[(1-e^{-i\xi})^\alpha+(1-e^{i\xi})^\alpha\Big],\cr
  M_o(\xi)&=    \frac{1}{2i\sin (\alpha \pi/2)\sgn(\xi)}
    \Big[(1-e^{-i\xi})^\alpha-(1-e^{i\xi})^\alpha\Big].
\end{align*}
By choosing the principal branch of the power function $(1-e^{\pm
i\xi})^{\alpha}$, we can further verify that $M_e(\xi) \sim |\xi|^{\alpha}$
and $M_o(\xi) \sim |\xi|^{\alpha}$
as $|\xi| \to 0$, and confirm that  the scheme~\eqref{eq:fd} with
the weights $\{w_k^{GL}\}$
provides a reasonable discretization of $\mathcal{A}$.

The case when $\alpha \in (1,2)$ can be worked out similarly,
starting from the fractional derivatives
\[
  \mathcal{D}_\pm^\alpha f(x)  = \frac{1}{\Gamma(2-\alpha)}\frac{\de^2}{\de x^2}
  \int_0^\infty t^{1-\alpha} f(x\mp t)\dt
  =\frac{1}{\Gamma(-\alpha)(2^{\alpha}-2)}\int_0^\infty \frac{f(x)-2f(x\mp t)+f(x\mp
  2t)}{t^{1+\alpha}}\dt
\]
and their finite difference approximations
\[
  \mathcal{D}_{\pm}^\alpha f(x_j) \approx
  \mathcal{D}_{h\pm}^\alpha f_j
  :=h^{-\alpha}\sum_{k=0}^\infty (-1)^k\binom{\alpha}{k}f_{j\mp (k-1)},
\]
where the indices are shifted by one to preserve
the non-negativity of the coefficient of $f_{j\pm k}$ with $k\neq 0$.
This modification is important in time-dependent problems like random
walk models~\cite{gorenflo1998random}. Hence the weights collected
from~\eqref{eq:fraclevy} become
\[
  h^{\alpha}w^{GL}_j =
  \begin{cases}
    \frac{1+\beta}{2\cos (\alpha\pi/2)} (-1)^{j} \binom{\alpha}{j+1},\qquad &j>1, \cr
    -\frac{1}{2\cos(\alpha\pi/2)}\left[
    1-\beta+\frac{\alpha(\alpha-1)}{2}(1+\beta)
    \right], & j=1, \cr
    \frac{\alpha}{\cos (\alpha\pi/2)}, & j=0,\cr
        -\frac{1}{2\cos(\alpha\pi/2)}\left[
    1+\beta+\frac{\alpha(\alpha-1)}{2}(1-\beta)
    \right], & j=-1, \cr
    \frac{1-\beta}{2\cos (\alpha\pi/2)} (-1)^{j} \binom{\alpha}{-j+1},\qquad &j<-1.
  \end{cases}
\]
Similarly, the rescaled multipliers can also be obtained as
\begin{align*}
M_e(\xi) &=   \frac{1}{2\cos (\alpha\pi/2)}
    \Big[e^{i\xi}(1-e^{-i\xi})^\alpha    +e^{-i\xi}(1-e^{i\xi})^\alpha\Big], \cr
M_o(\xi) &=   \frac{1}{2i\sin (\alpha\pi/2)\sgn(\xi)}
    \Big[e^{i\xi}(1-e^{-i\xi})^\alpha    +e^{-i\xi}(1-e^{i\xi})^\alpha\Big].
\end{align*}
One reason for the popularity of this scheme is the non-negativity
of the weights for all $\beta \in [-1,1]$, important to
non-negativity of the probability densities and other stability
properties when used for evolutionary equations. However, the weights
becomes singular as $\alpha$ goes to $1$ (in either direction),
which might be suggested by above two distinct forms of $\mathcal{D}_\pm^\alpha$;
the case $\alpha = 1$ has to be treated differently,
for instance in~\cite{Gorenflo2002521}.

\subsection{Regularized Spectral weights}
Finally, we consider the scheme with
$M_e(\xi)=M_o(\xi)=\big(2-2\cos(\xi)\big)^{\alpha/2}$,
which is motivated from the rescaled multiplier $2-2\cos(\xi)$ for the
standard three-point central difference method for the second order
derivative and can be
thought as a regularised function of the exact multiplier $|\xi|^{\alpha}$ near the origin.
The scheme for the symmetric case ($\beta=0$) was already studied in~\cite{MR2365531},
together with different boundary conditions associated with a
bounded domain. The weights for general $\beta \in (-1,1)$  are
\begin{align*}
    w_j^{RS}
    &= -\frac{h^{-\a}}{\pi} \left[
      \int_0^\pi \big(2-2\cos(\xi)\big)^{\a/2} \cos(j\xi) \dxi
    -\beta \tan\left(\frac{\pi \a}{2}\right)
    \int_0^\pi \big(2-2\cos(\xi)\big)^{\a/2} \sin(j\xi) \dxi \right] \cr
      &= (1+\b)h^{-\a}\sin\left({\frac{\pi\a}{2}}\right)
      \frac{\Gamma(j-\frac{\a}{2})\Gamma(1+\a)}
      {\pi\Gamma(j+1+\frac{\a}{2})} \cr
&\quad      - \frac{(-1)^jh^{-\a}}{\pi (j-\frac{\a}{2})}\beta
      \tan\left(\frac{\pi \a}{2}\right)
      {}_2F_1\left(-\a,j-\frac{\a}{2};j+1-\frac{\a}{2};-1\right),
    \end{align*}
where detailed derivation  is given in Appendix~\ref{sec:wrs}.

For $\a= 1$, we have (by taking the limit as $\alpha$ goes to $1$)
    \[
      \int_0^\pi \sqrt{2-2\cos\xi}\cos (j\xi)\dxi =
      \frac{\Gamma(j-1/2)}{\Gamma(j+1+1/2)},
    \]
    but the other integral
    $\int_0^\pi \sqrt{2-2\cos\xi}\ln \xi \sin (k\xi)\dxi$
    does not seem to possess a closed form expression
    and hence has to be calculated using numerical quadrature~\cite{MR2211043}.

Here we only focus on schemes whose weights can be expressed explicitly as above,
while in principle more general ones can be constructed from given
pairs rescaled multipliers $M_{e}(\xi)$ and $M_{o}(\xi)$ that behaves like $|\xi|^{\alpha}$
near the origin. When the schemes are used in practical problems, other issues
like accuracy and non-negativity of the weights start to play important roles,
as commented in the numerical experiments in the next section.


\section{Numerical experiments}

In this section, several numerical experiments are performed to compare the proposed schemes,
for their accuracy and other related issued in the discretization of
practical problems.

\subsection{Accuracy of the schemes}
Because of the non-locality of the operator, the accuracy of different schemes is better
to be assessed through the spectral space.
Let $\{u_j\}$ be sampled from a function $u$ defined on $\mathbb{R}$.
If the Fourier transform  of $u$ is defined as
$\hat{u}(\xi) = \int_{-\infty}^\infty u(x)e^{-ix\xi}\mathrm{d}x$,
then
\[
  \mathcal{A}u(x_j) = \frac{1}{2\pi}\int_{-\infty}^\infty
  \psi(-\xi) \hat{u}(\xi)e^{ix_j\xi}\dxi,
  \quad \mathcal{A}_hu_j = \sum_{k=-\infty}^\infty w_{j-k}u_k
  =\frac{1}{2\pi}\int_{-\pi/h}^{\pi/h}
  \psi_h(-\xi)\hat{u}(\xi)e^{ix_j\xi}\dxi.
\]
Hence the error between $\mathcal{A}u(x_j)$ and its difference approximation
at $x_j$, provided that the infinite sum $\sum_k w_{j-k}u_k$ can be performed accurately, is given by
\begin{equation}\label{eq:errterm}
  \frac{1}{2\pi}\int_{|\xi|>\pi/h}   \psi(-\xi) \hat{u}(\xi)e^{ix_j\xi}\dxi
  + \frac{1}{2\pi}\int_{-\pi/h}^{\pi/h}
  \big( \psi(-\xi) - \psi_h(-\xi)\big)\hat{u}(\xi)e^{ix_j\xi}\dxi.
\end{equation}
If $\hat{u}(\xi)$ decays to zero fast enough as $|\xi|$ goes to infinity,
the error is dominated by the second integral. In other words,
the accuracy depends on how well $\psi(-\xi)$ is approximated by $\psi_h(-\xi)$,
or equivalently how well $|\xi|^\alpha$ is approximated by $M_e(\xi)$ and $M_o(\xi)$.

For the scheme with spectral weights $\{w_j^{SP}\}$, $M_e(\xi) = M_o(\xi) = |\xi|^\alpha$ and the
accuracy is expected to be spectral---the error usually decays like $O(e^{-C/h})$
for some positive constant $C$.

For the scheme with \GL weights $\{w_j^{GL}\}$, the rescaled multipliers
can be rewritten to better illustrate their behaviours near the the origin: for $\alpha \in (0,1)$,
\[
  M_e(\xi) = \frac{\cos \frac{\pi-|\xi|}{2}\alpha}{\cos \frac{\alpha\pi}{2}}\left(
  2\sin \frac{|\xi|}{2}\right)^{\alpha},\qquad
  M_o(\xi) = \frac{\sin \frac{\pi-|\xi|}{2}\alpha}{\sin \frac{\alpha\pi}{2}}\left(
  2\sin \frac{|\xi|}{2}\right)^{\alpha}
\]
and for $\alpha \in (1,2)$,
\[
  M_e(\xi) = \frac{\cos\left( \frac{\pi-|\xi|}{2}\alpha+|\xi|\right)}{\cos \frac{\alpha\pi}{2}}\left(
  2\sin \frac{|\xi|}{2}\right)^{\alpha},\quad
  M_o(\xi) = \frac{\sin\left( \frac{\pi-|\xi|}{2}\alpha+|\xi|\right)}{\sin \frac{\alpha\pi}{2}}\left(
  2\sin \frac{|\xi|}{2}\right)^{\alpha}.
\]
After some complicated algebra, all these multipliers can be expanded as
$|\xi|^{\alpha}(1+c_1\xi+c_2\xi^2+\cdots)$ for some non-zero constant $c_1$ (different
for different multipliers).
Therefore, the second integral in~\eqref{eq:errterm} becomes
\begin{multline*}
  \frac{1}{2\pi}\int_{-\pi/h}^{\pi/h}
    \big( \psi(-\xi) - \psi_h(-\xi)\big)\hat{u}(\xi)e^{ix\xi}\dxi
    = \frac{1}{2\pi}\int_{-\pi/h}^{\pi/h} \big( h^{-\alpha} M_e(h\xi)-|\xi|^{\alpha}\big)
    \hat{u}(\xi) e^{ix_j\xi}\dxi \cr
    + \frac{i\beta}{2\pi}\tan\left(\frac{\alpha\pi}{2}\right)
    \int_{-\pi/h}^{\pi/h} \big( h^{-\alpha} M_o(h\xi)-|\xi|^{\alpha}\big)
        \hat{u}(\xi)\sgn(\xi)e^{ix_j\xi}\dxi.
\end{multline*}
Upon substitution of the expansion for $M_e(\xi)$ and $M_o(\xi)$, the above error
is bounded by $O(h)$, provided that $\int_{-\pi/h}^{\pi/h} |\xi|^{\alpha+1}
|\hat{u}(\xi)|\dxi$ is finite.

Finally for the scheme with regularised spectral weights $\{w_j^{RS}\}$,
\[
  M_e(\xi) = M_o(\xi)= \big(2-2\cos\xi\big)^{\alpha/2} = |\xi|^{\alpha}\left(
  1 - \frac{\alpha}{24}\xi^2+\cdots
  \right).
\]
A similar procedure shows the leading $O(h^2)$ error,
provided that $\int_{-\pi/h}^{\pi/h} |\xi|^{\alpha+2}
|\hat{u}(\xi)|\dxi$ is finite.

\begin{figure}[htp]
\begin{center}
\includegraphics[totalheight=0.25\textheight]{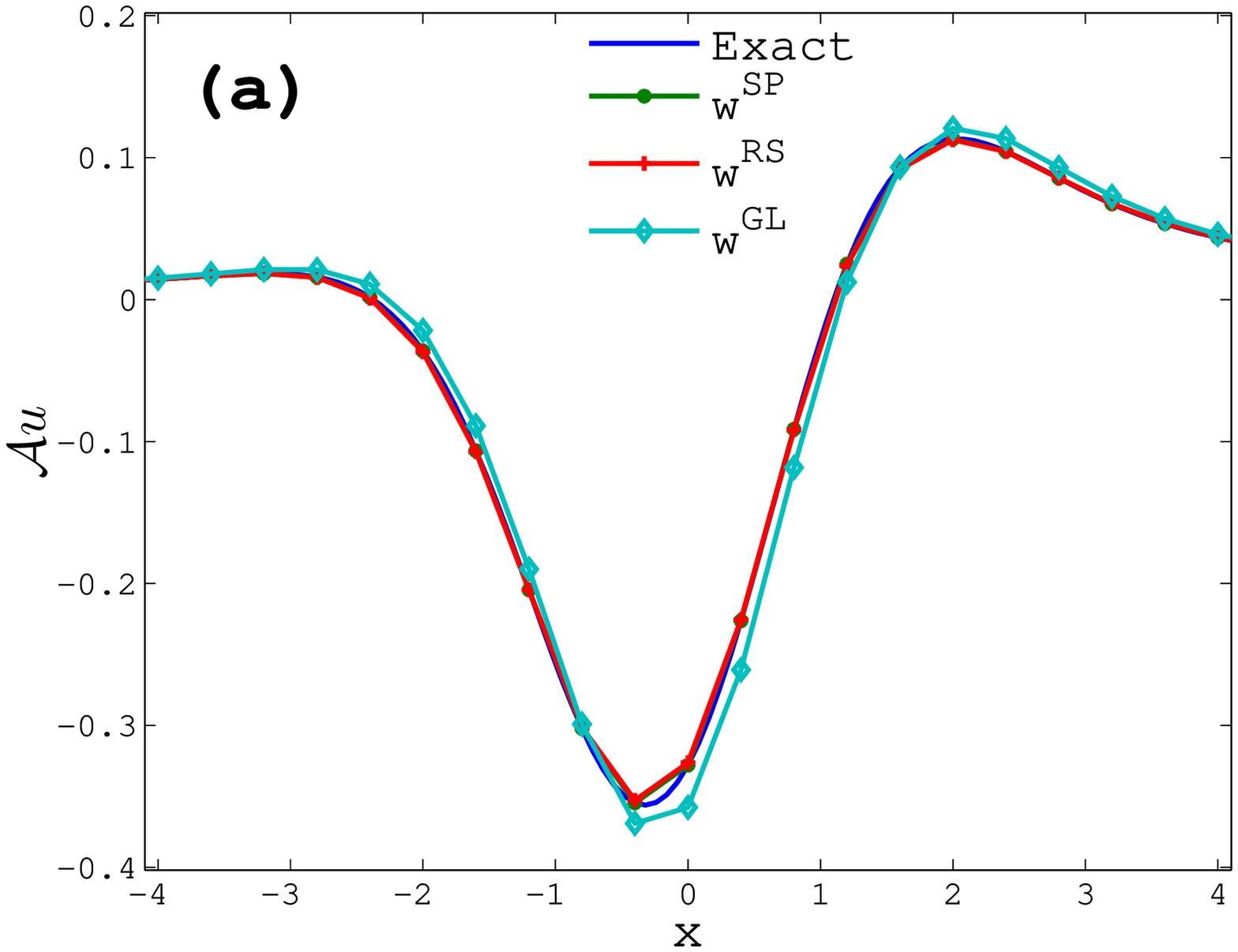}
$~ ~ ~ ~$
\includegraphics[totalheight=0.25\textheight]{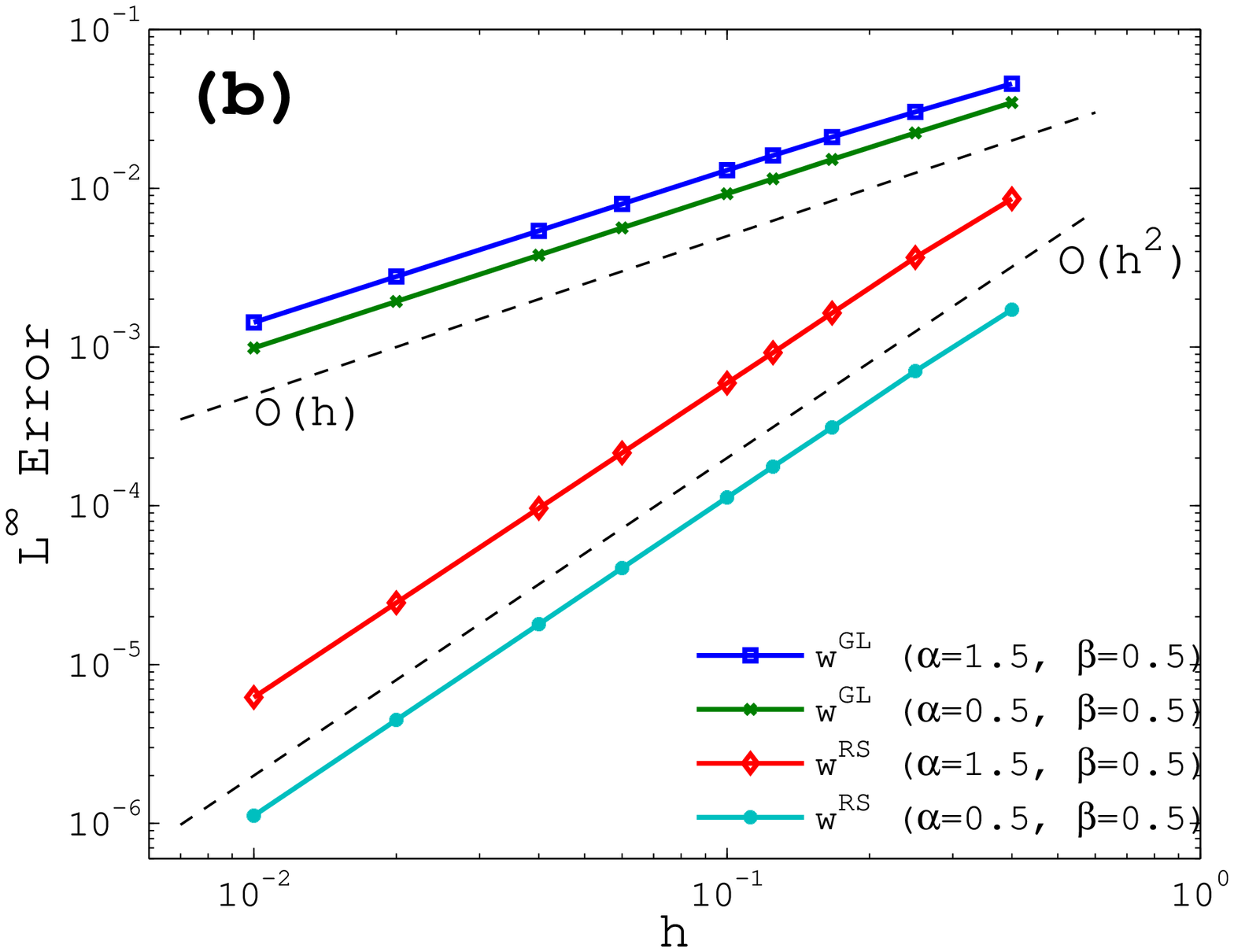}
\end{center}
\caption{ {\bf (a)} Different approximations of $\mathcal{A}u$ for
$u(x)=\frac{1}{\sqrt{2\pi}}e^{-x^2/2}$, and $\alpha=0.5,\beta=0.5$ with grid size $h=0.4$; {\bf (b)}
the convergence of the schemes with \GL weights and regularised spectral weights of the
same approximations with different grid size $h$.}
\end{figure}

Different approximations to $\mathcal{A}u$ (with $\alpha=0.5,\beta=0.5$)
for the Gaussian $u(x) = \frac{1}{\sqrt{2\pi}}e^{-x^2/2}$ on
a grid size $h=0.4$ are shown in Figure 1(a), where the 'exact'
approximation is computed using high order numerical quadrature for the inverse Fourier transform.
While small error could still be observed when the \GL weights are used,
the approximations are almost indistinguishable from the exact value $\mathcal{A}u$
for the spectral weights and regularised spectral weights,
even on a coarse grid with $h=0.4$.
The above analysis about the order of accuracy is further verified in Figure 1(b),
when the $L^\infty$ error (computed on the interval $[-4,4]$) decreases
with the expected order as the grid size $h$ is refined. The spectral convergence for $w^{SP}$
is omitted in the figure, because the error is already less than $10^{-6}$ on a grid size $h=0.4$.

Notice that the above convergence rates are only expected under certain restrictions:
the function $u(x)$ should be smooth enough such that
the integral $\int_{|\xi|>\pi/h} \psi(-\xi)\hat{u}(\xi)e^{ix\xi}\dxi$
can be safely ignored;  care must be taken to avoid introducing any error
when the infinite sum in the convolution of the scheme~\eqref{eq:fd}
is truncated.

\subsection{Application to mean exit time}

In practice, we are more interested in solutions to non-local differential equations
with the generator $\mathcal{A}$, rather than approximations
of the operator alone.
For example,  the mean exit time appears in many systems driven by various random
noises~\cite{completely_asy1,finite,first}.
Let the first exit time starting at $x$ from a bounded domain $\Omega$ is defined
as $ \tau{(\omega,x)}:= \inf \{t \geq 0,X_0=x, X_{t}(\omega , x) \notin  \Omega \}$, and the
mean first exit time (in short, mean exit time) is $u(x)=\mathbb{E}[\tau(\omega,x)].$
If $X_t$ is an $\alpha$-stable L\'evy motion, such that $\mathbb{E}\big(\exp(i\xi (X_t-X_0))\big)
=\exp(t\psi(\xi))$ with
\[
  \psi(\xi) = \begin{cases} -|\xi|^{\alpha}\big( 1-i\beta \sgn(\xi)\tan(\alpha\pi/2)\big),\qquad
&    \alpha \neq 1,\cr
    -|\xi|\big(1-i2\beta \sgn(\xi) \ln |\xi|/\pi\big), & \alpha = 1,
  \end{cases}
\]
then the mean exit time satisfies the following nonlocal partial differential equation \cite{Duan}
\begin{align}\label{eq:exit2D}
  \mathcal{A} u(x) &= -1, \quad \text{for} \; x \in \Omega,
\end{align}
  subject to the Dirichlet-type exterior condition $u(x) \equiv 0$ for $x \in \Omega^c$.

Because of the zero boundary condition outside the domain, $u_j \equiv 0$ for $|j|\geq 1/h$ and
the unknowns $\{u_j\}$ are governed by a linear system of equations
\[
  \sum_{|k|<1/h} w_{j-k}u_k = -1,\quad |j|<1/h.
\]
This system is Toeplitz and can usually be solved very efficiently~\cite{MR2108963}.

\begin{figure}[htp]
\begin{center}
\includegraphics[totalheight=0.25\textheight]{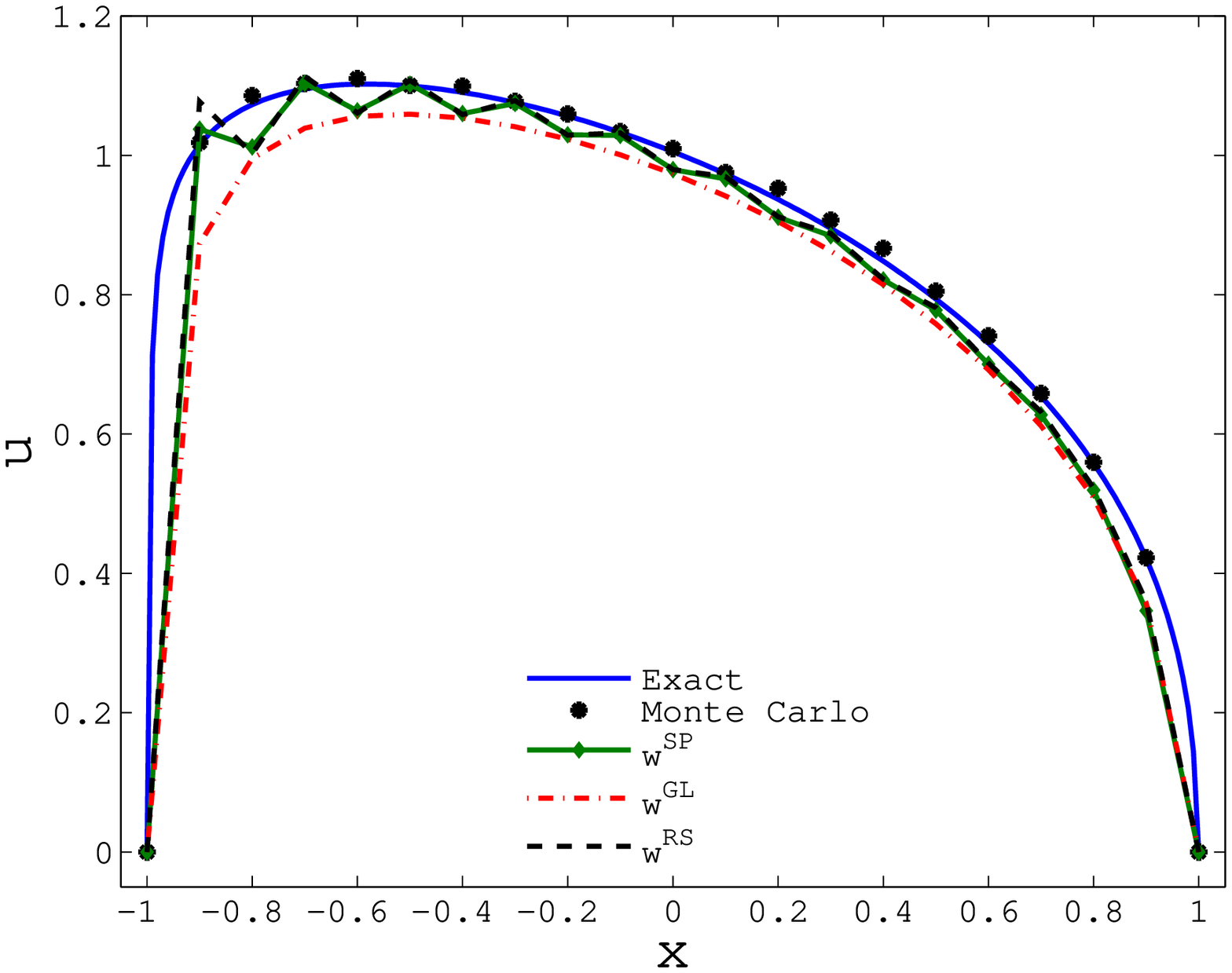}
$~ ~ ~ ~$
\includegraphics[totalheight=0.25\textheight]{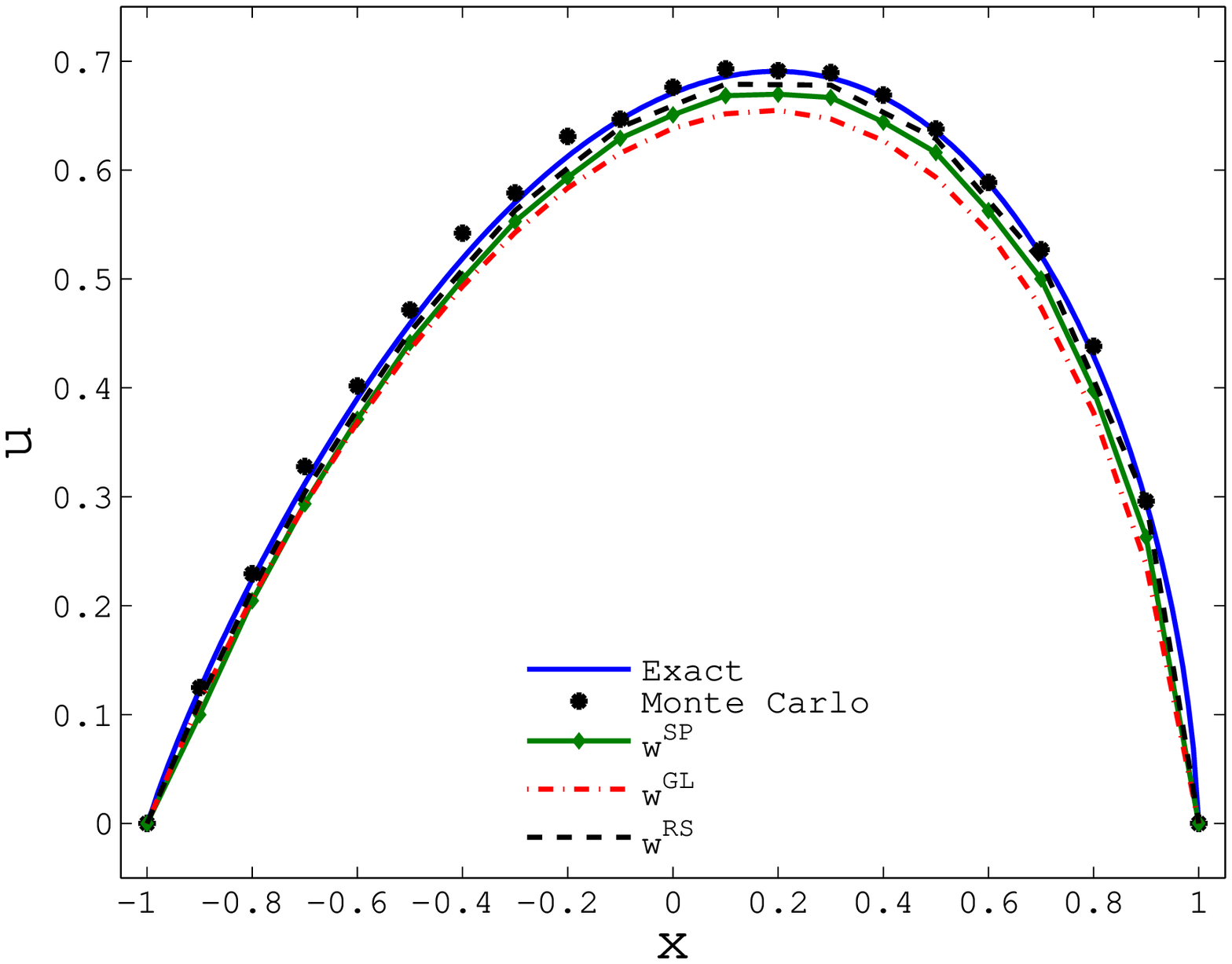}
\end{center}
\caption{The solution to the mean exit problem on $\Omega=[-1,1]$, where the results
  from different schemes are compared  with the time estimated from Monte Carlo simulations:
  $\alpha=0.5,\beta=0.5$ (left
figure) and $\alpha = 1.5,\beta=0.5$(right figure). The grid size $h$ is $0.1$ and
the 'exact' solution is obtained from the numerical solution on a refined grid $h=0.001$
using \GL weights.}
\end{figure}

The numerical solutions for $\alpha = 0.5$ and $\alpha = 1.5$ ($\beta=0.5$ in both cases)
are show in Figure 2 with grid size $h=0.1$. As a comparison, mean exit time estimated from Monte Carlo
simulation (see~\cite{MR1306279} for related random number generation) is also plotted,
with 10000 sample paths for each process starting at $x$
and time step $0.0001$ for advancing the evolution. While all numerical solutions
agree with the Monte Carlo simulation to some degree,
new problems arise as a result of the non-smoothness of the
solution $u$ near the boundaries $x=\pm 1$. {For instance, the mean exit time $u(x)$
  for the symmetric operator behaves
  like $\mbox{dist}(x,\partial\Omega)^{\alpha/2}$ near the boundary~\cite{ROSOTON2014275},
which is not even Lipschitz continuous for any $\alpha \in (0,2)$.}
Because of this loss of regularity, the order of convergence derived in the
previous subsection for each scheme no longer holds.
Moreover, because of the appearance of negative weights,
the numerical solutions become oscillatory, which is more pronounced
when $\alpha$ is less than one, or when $|\beta|$ is close to one.

The situation is even worse for evolutionary {problems}
when schemes with negative weights (except $w_0$, which is always negative yet does not
affect the stability)
are used---probability densities~\footnote{Strictly speaking, it is the formal adjoint
  operator $\mathcal{A}^*$ should be used for the evolution of probability densities.}
  could become negative during the evolution. In cases the underlying process
  is genuinely asymmetric ($\beta\neq 0$), it is common to have negative weights,
because the magnitude of the integral $\int_0^\pi M_o(\xi)\sin(j\xi)\dxi$
usually decays to zero slower than that of  $\int_0^\pi M_e(\xi)\cos(j\xi)\dxi$.
In this regard, only the scheme with \GL weights has the desired stability property
for all $\beta \in [-1,1]$.

\section{Conclusion}

In this paper, we proposed several difference schemes of the general form~\eqref{eq:fd}
for the infinitesimal generator of the $\alpha$-stable L\'evy process,
mainly based on the corresponding representations under the semi-discrete Fourier
transform, such that the exact multiplier is properly approximated.
The scheme with spectral weights $\{w_j^{SP}\}$ is spectrally accurate for smooth functions, but the 
accuracy becomes degenerate for non-smooth functions, and the solutions
show spurious {oscillations} in the {discretization} of practical problems.
The scheme with classical \GL weights $\{w_j^{GL}\}$ exhibits good numerical stability, but
it is only first order accurate and a ``singularity'' appears when $\alpha$ approaches one.
A {compromise} between accuracy and stability
is reached in the scheme with regularised spectral weights
$\{w_j^{RS}\}$: the accuracy is second order, adequate for most problems, while
non-physical oscillations are not as pronounced as those when the spectral weights are used.
The construction of high order schemes with good stability stability remains
a challenging task. {In practice, if the function under consideration
  does not decay to zero fast enough,
the accuracy of the scheme could be reduced when the discrete operator $\mathcal{A}_h$
is  just truncated to finite sums.
The next order may still be recovered in many circumstances, by combining far-field
asymptotic behaviour of the function and the integral representation in~\eqref{eq:levyinfgen},
in a similar way as in~\cite[Section 5]{HuangOberman1}. However, this treatment of
boundary condition requires certain information like the decay rate of the function
that may demand deep theoretical analysis.
}

\appendix

\section{The generalised hypergeometric function}
\label{sec:spint}

The generalised hypergeometric function appears at several places in this paper,
and near the origin, it is defined as a series
\begin{equation}\label{eq:ghg}
  { }_pF_q\big(a_1,a_2,\cdots,a_p; b_1,b_2,\cdots,b_q; z\big)
  =\sum_{n=0}^\infty \frac{(a_1)_n(a_2)_n\cdots (a_p)_n}{
  (b_1)_n(b_2)_n\cdots (b_q)_n n!} z^n,
\end{equation}
where $(a)_n = a(a+1)\cdots (a+n-1)$ is the Pochhammer symbol.
The most common ones include Kummer's confluent hypergeometric
function ${}_1F_1$ and Gauss's hypergeometric function ${ }_2F_1$.
From the series representation~\eqref{eq:ghg}, it is easy to verify that
\begin{multline}\label{eq:ghgrec}
  \frac{\mathrm{d}}{\mathrm{d}z}  { }_pF_q\big(a_1,a_2,\cdots,a_p; b_1,b_2,\cdots,b_q; z\big) \cr
= \frac{a_1a_2\cdots a_p}{b_1b_2\cdots b_q}
  { }_pF_q\big(a_1+1,a_2+1,\cdots,a_p+1; b_1+1,b_2+1,\cdots,b_q+1; z\big).
\end{multline}

\section{Special definite integrals} \label{sec:spin}
In the simplification of the function $\psi(\xi)$ when the L\'evy measure
is given by~\eqref{eq:alphalevymeas}, the
following special definite integrals are used (see \cite[p.436-437]{MR1398882}):
\begin{subequations}
  \label{eq:intas}
\begin{align}
  \int_0^\infty x^{-\alpha}\sin x\ \mathrm{d}x &= \Gamma(1-\alpha)\cos
  \left(\frac{\pi \alpha}{2}\right),  \qquad 0<\alpha < 2\\
  \int_0^\infty x^{-\alpha}\cos x\ \mathrm{d}x &= \Gamma(1-\alpha)\sin
  \left(\frac{\pi \alpha}{2}\right),  \qquad 0<\alpha<1,
\end{align}
together with other variants like
\begin{equation}
  \int_0^\infty x^{-\alpha}(1-\cos x) \mathrm{d}x = -\Gamma(1-\alpha)\sin
  \left(\frac{\pi \alpha}{2}\right),  \qquad 1<\alpha<3,
\end{equation}
and
\begin{equation}
  \int_0^1 \frac{\sin x-x}{x^2}\mathrm{d}x + \int_1^\infty \frac{\sin x}{x^2} \mathrm{ d}x= 1-\gamma,
\end{equation}
where $\gamma$ is the Euler-Mascheroni constant.
\end{subequations}

\section{Evaluation of the integrals in the weights $w^{RS}_j$}
\label{sec:wrs}
Here  explicit expressions of the weights $\{w_j^{RS}\}$
associated with the regularised
multipliers $M_e(\xi) =M_o(\xi)= \big(2-2\cos \xi)^{\alpha/2}$
are derived. In fact, $M_e(\xi)=M_o(\xi)= 2^{\alpha}\sin^{\alpha}\left(\frac{\xi}{2}\right)$,
and   we can verify the following indefinite integral
\begin{equation}\label{eq:rsqint}
  \int \sin^{\alpha}\left( \frac{\xi}{2}\right)e^{ij\xi} \de \xi
  =\frac{i^{\alpha+1}}{2^{\alpha}(\alpha/2-j)}
  e^{i(j-\alpha/2)\xi}{ }_2F_1\left(-\alpha,j-\frac{\alpha}{2};
  j+1-\frac{\alpha}{2}; e^{i\xi} \right).
\end{equation}
First, using~\eqref{eq:ghgrec} for the derivative
of generalised hypergeometric functions,  we get
\begin{multline} \label{eq:dgf}
  \  \frac{e^{-i(j-\alpha/2)\xi}}{i(j-\alpha/2)}\frac{\de}{\de\xi}
  e^{i(j-\alpha/2)\xi}
  { }_2F_1\left(-\alpha,j-\frac{\alpha}{2}; j+1-\frac{\alpha}{2}; e^{i\xi} \right) =\cr
      { }_2F_1\left(-\alpha,j-\frac{\alpha}{2}; j+1-\frac{\alpha}{2}; e^{i\xi} \right)
      -\frac{\alpha e^{i\xi}}{j+1-\alpha/2}
  { }_2F_1\left(1-\alpha,j+1-\frac{\alpha}{2}; j+2-\frac{\alpha}{2}; e^{i\xi} \right).
\end{multline}
By the series representation of ${ }_2F_1$,
the right hand side of~\eqref{eq:dgf} becomes
\begin{align*}
    &\quad 1 +(j-\alpha/2)\sum_{n=1}^\infty \frac{(-\alpha)_n}{n!(j+n-\alpha/2)}e^{in\xi}
      -\alpha\sum_{n=0}^\infty \frac{(1-\alpha)_n}{n!(j+1+n-\alpha/2)}
      e^{i(n+1)\xi} \cr
    &= 1 +(j-\alpha/2)\sum_{n=1}^\infty \frac{(-\alpha)_n}{n!(j+n-\alpha/2)}e^{in\xi}
    -\alpha\sum_{n=1}^\infty \frac{(1-\alpha)_{n-1}}{(n-1)!(j+n-\alpha/2)}
      e^{in\xi} \cr
    &= 1 + \sum_{n=1}^\infty \frac{(-\alpha)_n}{n!}e^{in\xi},
\end{align*}
which is exactly $(1-e^{i\xi})^\alpha$. Therefore, the
indefinite integral~\eqref{eq:rsqint} is established by
\begin{multline*}
  \quad \frac{\de }{\de \xi}\left[  \frac{i^{\alpha+1}}{2^{\alpha}(\alpha/2-j)}
    e^{i(j-\alpha/2)\xi}{ }_2F_1\left(-\alpha,j-\frac{\alpha}{2};
  j+1-\frac{\alpha}{2}; e^{i\xi} \right)\right] \cr
  = \left(\frac{i}{2}\right)^\alpha e^{-i\alpha \xi/2}(1-e^{i\xi})^{\alpha}e^{ij\xi}
  = \left[ \frac{e^{i\xi/2}-e^{-i\xi/2}}{2i}\right]^{\alpha}e^{ij\xi}
  = e^{ij\xi}\sin^\alpha\left( \frac{\xi}{2}\right) ,
\end{multline*}
while in the last step the principal branch of the fractional power
is taken with $\xi$ assumed to be in the interval $[0,\pi]$.

Extract the real and the imaginary parts of the following integral
\begin{multline*}
  \int_0^{\pi} \big(2-2\cos\xi\big)^{\alpha/2} e^{ij\xi}\dxi
  =\frac{i^{\alpha+1}}{\alpha/2-j}\left[
    e^{i(j-\alpha/2)\pi}
  { }_2F_1\left(-\alpha,j-\frac{\alpha}{2};j+1-\frac{\alpha}{2};-1
  \right)\right. \cr
  \left. -{
  }_2F_1\left(-\alpha,j-\frac{\alpha}{2};j+1-\frac{\alpha}{2};1\right)
\right],
\end{multline*}
we get
\begin{equation*}\label{eq:rscos}
  \int_0^\pi \big(2-2\cos\xi\big)^{\alpha/2} \cos (j\xi) \de \xi
  =-\frac{\Gamma(j-\alpha/2)\Gamma(1+\alpha)}{\Gamma(j+1+\alpha/2)}
  \sin \frac{\alpha\pi}{2},
\end{equation*}
and
\begin{multline*}\label{eq:rssin}
  \int_0^\pi \big(2-2\cos\xi\big)^{\alpha/2} \sin (j\xi) \de \xi
  =\frac{\Gamma(j-\alpha/2)\Gamma(1+\alpha)}{\Gamma(j+1+\alpha/2)}\cos
  \frac{\alpha\pi}{2} \cr
  +\frac{(-1)^j}{\alpha/2-j} { }_2F_1\left(-\alpha,j-\frac{\alpha}{2};
  j+1-\frac{\alpha}{2};-1\right),
\end{multline*}
where are simplified by Gauss's identity
\[
  { }_2F_1\left(-\alpha,j-\frac{\alpha}{2};
  j+1-\frac{\alpha}{2};1\right) =
  \frac{\Gamma(j+1-\alpha/2)\Gamma(1+\alpha)}{\Gamma(j+1+\alpha/2)}.
\]

\end{document}